\documentclass [12pt,leqno] {article}
\usepackage{amsmath,amsfonts,url}
\usepackage{enumerate}
\usepackage{color}

\newcommand{\normal}{\color{black}}

\def\supp{\mathop{\rm supp}\nolimits}

\newcommand{\dowod}{{\em Proof}.\/ }
\newcommand{\dowodof}{{\em Proof}}
\newcommand{\qed}{~\hfill~$\fbox{}$ \vspace{0.5cm}}
\newcommand{\R}{ \mathbb{R}^{d}}
\newcommand{\N}{\mathbb{N}}
\newcommand{\indyk}[1]{{\bf 1}_{#1}}

\newcommand{\gener}{{\cal A}}
\newcommand{\scalp}[2]{#1\cdot#2}

\newcommand{\supmasyL}{\bar{b}_\varepsilon}

\newcommand{\infmasyL}{\underline{b}_\varepsilon}
\newcommand{\Cvanish}{C_\infty(\R)}
\newcommand{\Bbound}{B_b(\R)}
\newcommand{\Aapprox}[1]{{\cal A}_{#1}}

\newtheorem{lemma}{\indent\sc Lemma}
\newtheorem{proposition}{\indent\sc Proposition}
\newtheorem{theorem}{\indent\sc Theorem}
\newtheorem{remark}{\indent\sc Remark}

\begin{document}

\title{Upper estimates of transition densities for stable-dominated semigroups }
\author{Kamil Kaleta and Pawe{\l} Sztonyk}
\footnotetext{ Kamil Kaleta \\ Institute of Mathematics,
  University of Warsaw, \\
  ul. Banacha 2,
  02-097 Warszawa, Poland.\\
}
\footnotetext{ Kamil Kaleta and Pawe{\l} Sztonyk \\ Institute of Mathematics and Computer Science,
  Wroc{\l}aw University of Technology,
  Wybrze{\.z}e Wyspia{\'n}\-skie\-go 27,
  50-370 Wroc{\l}aw, Poland.\\
  {\rm e-mail: Kamil.Kaleta@pwr.wroc.pl, Pawel.Sztonyk@pwr.wroc.pl} \\
}
\maketitle

\begin{abstract}
  We derive upper estimates of transition densities for Feller semigroups with jump intensities lighter
  than that of the rotation invariant stable L\'evy process.  
\end{abstract}
\footnotetext{2000 {\it MS Classification}:
Primary 60J75, 60J35; Secondary 47D03.
}

\section{Introduction and Preliminaries}

Let $\alpha\in(0,2)$ and $d=1,2,\dots$. For the rotation invariant $\alpha$ - stable L\'evy process on $\R$ with the L\'evy measure
\begin{equation}\label{eq:RISLM}
  \nu(dy)=\frac{c}{|y|^{\alpha+d}}\,dy,\quad y\in\R\setminus\{0\},
\end{equation}
the asymptotic behaviour of its transition densities $p(t,x,y)$ is well-know (see, e.g., \cite{BlGe}), i.e.
$$
  p(t,x,y) \approx \min\left(t^{-d/\alpha},\frac{t}{|y-x|^{\alpha+d}}\right),\quad t>0,\,x,y\in\R.
$$

Estimates of densities for more general classes of stable and other jump L\'evy processes gradually extended. Obtained results
contained estimates for general stable processes in \cite{W,BS2007} and tempered and layered stable processes in \cite{Szt2008} 
and \cite{Szt2011}.

In \cite{Szt2010} estimates of semigroups of stable-dominated Feller operators are given. The corresponding Markov process
is a Feller process and not necessarily a L\'evy process.
The name {\it stable dominated} refers to
the fact that the intensity of jumps for the investigated semigroup is dominated by (\ref{eq:RISLM}). In the present paper
we extend the results obtained in \cite{Szt2010} and give estimates from above for a wider class of semigroups with intensity of jumps
lighter than stable processes. We will now describe our results.

Let $f:\R\times\R\mapsto [0,\infty]$ be a Borel function. We consider the following assumptions on $f$. 

\bigskip
\noindent
\textbf{(A.1)} There exists a constant $M>0$ such that
$$
  f(x,y) \leq M \frac{\phi(|y-x|)}{|y-x|^{\alpha+d}} ,\quad x,y\in\R,\, y\neq x,
$$
where $\phi:\: [0,\infty)\to (0,1]$ is a Borel measurable function such that 
\begin{itemize}
\item[(a)] $\phi(a)=1$ for $a \in [0,1]$ and there is a constant $c_1=c_1(\phi)$ such that
$$
  \phi(a) \leq c_1 \phi(b),\quad |a-b| \leq 1 ,
$$
\item[(b)] $\phi \in C^2(1,\infty)$ and there is a constant $c_2=c_2(\phi,\alpha,d)$ such that
$$
 \max\left(\left|\phi'(a)\right|,\left|\phi''(a)\right|\right) \leq c_2 \phi(a) 
$$
for every $a > 1$. 
\item[(c)] there is $c_3=c_3(\phi,\alpha,d)$ such that
$$
  \int_{|x-z|\geq 1, |y-z|\geq 1} \frac{\phi(|y-z|)}{|y-z|^{\alpha+d}} \frac{\phi(|z-x|)}{|z-x|^{\alpha+d}}\,dz \leq c_3 \frac{\phi(|y-x|)}{|y-x|^{\alpha+d}},
$$
for every $|x-y|>2$.
\end{itemize}

\bigskip
\noindent
\textbf{(A.2)}  $\,f(x,x+h)=f(x,x-h)\,$ for all $x,h\in\R$, or $\alpha < 1$.

\bigskip
\noindent
\textbf{(A.3)} $f(x,y)=f(y,x)$ for all $x,y\in\R$.

\bigskip
\noindent
\textbf{(A.4)} There exists a constant $c_4=c_4(\phi,\alpha,d)$ such that 
$$
  \inf_{x\in\R} \int_{|y-x|>\varepsilon} \frac{f(x,y)}{\phi(|y-x|)}\,dy \geq c_4 \varepsilon^{-\alpha},\quad \varepsilon >0.
$$
\noindent Denote 
$$
  b_\varepsilon(x)=\int_{|y-x|>\varepsilon} f(x,y)\,dy,\,\quad \varepsilon>0,\,x\in\R.
$$
It follows from (A.1) that there is also the constant $c_5=c_5(\phi,\alpha,d)$ such that 
$$
  \supmasyL := \sup_{x\in\R}b_\varepsilon(x) \leq c_5 \varepsilon^{-\alpha},\quad 0 < \varepsilon \leq 1.
$$
Thus, (A.4) is a partial converse of (A.1) and we have
$$
  \infmasyL := \inf_{x\in\R}b_\varepsilon(x) \geq c_6 \varepsilon^{-\alpha},\quad 0 < \varepsilon \leq \varepsilon_0\normal, 
$$
 for constants $\varepsilon_0=\varepsilon_0(\phi,\alpha,d), c_6=c_6(\phi,\alpha,d)$. 

We note that the assumption (A.1)(c) is satisfied for every nonincreasing function $\phi: (0,\infty)\to (0,1]$ such that
$$
  \phi(a)\phi(b) \leq c\, \phi(a+b),\quad a,b >1,
$$
for some positive constant $c$. Therefore it is easy to verify that all the assumptions on $\phi$ are satisfied, e.g.,
for functions $\phi(s)=e^{(1-s^{\beta})} \wedge 1$, where $\beta\in (0,1]$, $\phi(s) = (1\vee s)^{-\gamma}$, where $\gamma>0$, 
$\phi(s) = 1/\log(e(s\vee 1))$, $\phi(s)=1/\log\log (e^e(s\vee 1))$, and all their products and positive powers. 

It is also reasonable to ask if the conditions in the assumption (A.1) are satisfied by more general functions of the form
\begin{align}
\label{eq:specphiform}
\phi(s) =
\left\{
\begin{array}{ll} 1
& \mbox{if $s \in [0,1]$,}\\
e^{-ms^{\beta}} s^{\gamma}  & \mbox{if $s > 1$, with $m, \beta>0$, $\gamma \in \mathbb{R}$}.
\end{array} \right.
\end{align}
In this case, both conditions (a) and (b) on $\phi$ hold for $\beta \in (0,1]$ with no further restrictions on parameters $m$ and $\gamma$, while, as proven in Section \ref{sec:Examples}, the condition (c) is satisfied when $\beta \in (0,1]$ and $\gamma < d/2+\alpha-1/2$. Furthermore, this restriction on parameters is essential (see Remark \ref{rm:rm} in Section 3). Note also that this range of $\beta$ and $\gamma$ in \eqref{eq:specphiform} covers, e.g., jump intensities dominated by those of isotropic relativistic stable processes (see e.g. \cite[Lemma 2.3]{KulSiu}). 
 
For $x\in\R$ and $r>0$ we let $B(x,r)=\{y\in\R :\:|y-x|<r\}$. $B_b(\R)$ denotes the set of bounded Borel measurable functions, $C^k_c(\R)$ denotes the set of $k$ times continuously differentiable functions with compact support and $\Cvanish$ is the set of 
continuous functions vanishing at infinity. 
We use $c,C$ (with subscripts) to denote finite positive constants
which depend only on $\phi$ (the constant $M$), $\alpha$ and the dimension $d$. Any {\it additional} dependence
is explicitly indicated by writing, e.g., $c=c(n)$.
The value of $c,C$, when used without subscripts, may
change from place to place. We write $f(x)\approx g(x)$ to indicate that there is a constant $c$ such that $c^{-1}f(x) \leq g(x) \leq c f(x)$. 

Under the assumptions (A.1) and (A.2) we may consider the operator
\begin{eqnarray*}
\gener\varphi(x) & =  &
\lim_{\varepsilon\downarrow 0}\int_{|y-x|>\varepsilon} \left(\varphi(y)-\varphi(x)\right)f(x,y)\,dy\\
&  =   & \int_{\R} \left(\varphi(x+h)-\varphi(x)-\scalp{h}{\nabla\varphi(x)}\indyk{|h|<1}\right)
          f(x,x+h)\,dh \\
&      & +\, \frac{1}{2}\int_{|h|<1} \scalp{h}{\nabla\varphi(x)} \left(f(x,x+h)-f(x,x-h) \right)\,dh ,\quad \varphi\in C^2_c(\R).
\end{eqnarray*}

Recall the following basic fact (see \cite[Lemma 1]{Szt2010}).
\begin{lemma}
If (A.1), (A.2) hold and the function
$x\to f(x,y)$ is continuous on $\R\setminus\{y\}$ for every $y\in\R$ then $\gener$ maps $C^2_c(\R)$ into $\Cvanish$.
\end{lemma}

In the following we always assume that the condition (A.1) is satisfied.
For every $\varepsilon>0$ we denote
$$
  f_{\varepsilon}(x,y)=\indyk{B(0,\varepsilon)^c}(y-x)f(x,y),\quad x,y\in\R,
$$
and
$$
  \Aapprox{\varepsilon} \varphi(x)=\int \left(\varphi(y)-\varphi(x)\right)f_{\varepsilon}(x,y)\,dy,\quad \varphi\in \Bbound.
$$

Note that the operators $\Aapprox{\varepsilon}$ are bounded since 
$|\Aapprox{\varepsilon}\varphi(x)| \leq 2 \| \varphi \|_\infty b_\varepsilon(x)\leq 2\supmasyL \| \varphi \|_\infty$.
Therefore the operator 
$$
  e^{t\Aapprox{\varepsilon}}=\sum_{n=0}^\infty \frac{t^n\Aapprox{\varepsilon}^n}{n!},\quad t\geq 0,\,\varepsilon>0,
$$
is well--defined and bounded from $\Bbound$ to $\Bbound$. In fact for every $\varepsilon>0$ the family of operators $\{e^{t\Aapprox{\varepsilon}},\,t\leq 0\}$ 
is a semigroup on $\Bbound$, i.e., $e^{(t+s)\Aapprox{\varepsilon}}=e^{t\Aapprox{\varepsilon}}e^{s\Aapprox{\varepsilon}}$ 
for all $t,s\geq 0$, $\varphi\in\Bbound$. We note that $e^{t\Aapprox{\varepsilon}}$ is positive for all $t\geq 0$, $\varepsilon>0$ (see (\ref{eq:exp})).

Our first result is the following theorem.

\begin{theorem}\label{th:MainT} 
If (A.1) -- (A.4) are satisfied then there exist the constants $C_1$ and $C_2$ such that 
for every nonnegative $\varphi\in \Bbound$ and $\varepsilon\in (0,\varepsilon_0\wedge 1)$ we have
$$
  e^{t\Aapprox{\varepsilon}}\varphi(x) \leq C_1 e^{C_2t} \int 
  \varphi(y) \min\left(t^{-d/\alpha},\frac{t \phi(|y-x|)}{|y-x|^{\alpha+d}}\right) dy + e^{-tb_\varepsilon(x)}\varphi(x),
$$
for every $x\in\R.$
\end{theorem}

The proof of Theorem \ref{th:MainT} is given in Section \ref{Est}. To study a limiting semigroup 
we will need some additional assumptions.

\bigskip
\noindent
\textbf{(A.5)} The function
$x\to f(x,y)$ is continuous on $\R\setminus\{y\}$ for every $y\in\R$.

\bigskip
\noindent
\textbf{(A.6)} $\gener$ regarded as an operator on $\Cvanish$ is closable and its closure $\bar{\gener}$ is a generator of a strongly continuous 
contraction semigroup of operators $\{P_t,\,t\geq 0\}$ on $\Cvanish$.

\bigskip
Clearly, for every $\varphi\in C^2_c(\R)$ with $\sup_{x\in\R}\varphi(x)=\varphi(x_0)\geq 0$ 
we have $\gener\varphi(x_0)\leq 0$, i.e., $\gener$ satisfies {\it the positive maximum principle}. This implies that all $P_t$ ($t \geq 0$) are positive operators (see \cite[Theorems 1.2.12 and 4.2.2]{EK}). Thus, by our assumptions, $\{P_t,\,t\geq 0\}$ is {\it a Feller semigroup}.

The following theorem is our main result.

\begin{theorem}\label{th:main}
If (A.1)--(A.6) hold then there is $p:(0,\infty)\times\R\times\R\to [0,\infty)$
 such that
$$
  P_t \varphi (x) = \int_{\R} \varphi(y) p(t,x,y)\, dy,\quad x\in\R,\,t>0,\,\varphi\in \Cvanish,
$$
and 
\begin{equation}\label{eq:main}
  p(t,x,y) \leq C_1 e^{C_2t} \min\left(t^{-d/\alpha}, \frac{t \phi(|y-x|)}{|y-x|^{\alpha+d}}\right),\quad x,y\in\R,\,t>0.
\end{equation}
\end{theorem}

We note that $\gener$ is conservative, i.e., for $\varphi\in C^{\infty}_c(\R)$ such that $0\leq \varphi\leq 1$, $\varphi(0)=1$, and 
$\varphi_k(x)=\varphi(x/k)$, we have
$\sup_{k\in\N}\|\gener\varphi_k\|_\infty<\infty$, and 
$\lim_{k\to\infty}(\gener\varphi_k)(x)=0$, for every
$x\in\R$. It follows from Theorem 4.2.7 in \cite{EK} that 
there exists a Markov process $\{X_t,\,t\geq 0\}$ such that
${\mathbb E}[\varphi(X_t)|X_0=x]=P_t\varphi(x)$.

It is known that every generator $G$ of a Feller semigroup with $C_c^\infty(\R)\subset {\cal D}(G)$
is necessarily of the form
\begin{eqnarray}\label{eq:FellerGen}
  G\varphi(x)
  &  =  & \sum_{i,j=1}^d q_{ij}(x)D_{x_i}D_{x_j} \varphi(x) + l(x) \nabla \varphi(x) -c(x)\varphi(x) \\
  &     & + \int\limits_{\R}\left(\varphi(x+h)-\varphi(x)-
            \scalp{h}{\nabla \varphi(x)}\;\indyk{|h|<1}\right)\, \nu(x,dh)\,, \nonumber
\end{eqnarray}
where $\varphi\in C_c^\infty(\R)$, $q(x)=(q_{ij}(x))_{i,j=1}^n$ is a nonnegative definite real symmetric matrix,
the vector $l(x)=(l_i(x))_{i=1}^d$ has real coordinates, $c(x)\geq 0$,
and $\nu(x,\cdot)$ is a L\'evy measure (see \cite[Chapter 4.5]{Jc1}).

The converse problem whether a given operator $G$ generates a Feller semigroup is not completely resolved yet.
For the interested reader we remark that criteria are given, 
e.g., in \cite{Hoh1,Hoh2,Jacob1,Jacob2,Jacob3}. Generally, smoothness of the coefficients $q,l,c,\nu$
in (\ref{eq:FellerGen}) is sufficient for the existence (see Theorem 5.24 in \cite{HohHab}, Theorem 4.6.7 in \cite{Jacob4} and Lemma 2 in \cite{Szt2010}). Other conditions are given also in \cite{SchUe}.

Z.-Q.~Chen, P.~Kim and T.~Kumagai in \cite{ChKum,ChKum08,ChKimKum} investigate the case of symmetric jump--type Markov processes 
on metric measure spaces by using Dirichlet forms. Under the assumption that the corresponding jump
kernels are {\it comparable} with certain rotation invariant functions, they prove
the existence and obtain estimates of the densities (see Theorem~1.2 in \cite{ChKimKum}) analogous to (\ref{eq:main}). 
In the present paper we propose completely different approach which is based on general approximation scheme recently 
devised in \cite{Szt2010}. In Theorem \ref{th:main} we assume the estimate (A.1) from above but we use (A.4) as the only
estimate for the size of $f$ from below. We also emphasize that we obtain exactly $\phi(|x-y|)$ in (\ref{eq:main}) and 
from \cite{ChKum08,ChKimKum} follow estimates with $\phi(c|x-y|)$ for some constant $c\in (0,1)$. This seems to be essential especially in the case of exponentially localized L\'evy measures. Our general framework, including a layout of lemmas, is similar to that in \cite{Szt2010}. However, in the present case the decay of the jump intensity may be significantly lighter than stable and, therefore, much more subtle argument is needed. Note that the new condition (A.1)(c), which is pivotal for our further investigations, is necessary for the two-sided sharp bounds similar to the right hand side of (\ref{eq:main}). 

Other estimates of L\'evy and L\'evy-type transition densities are discussed in \cite{KnopKul,KnopSchill}.
In \cite{Lewand, Szt2007} the derivatives of stable densities have been considered, while bounds of heat kernels 
of the fractional Laplacian perturbed by gradient operators were studied in \cite{BJ2007}. An alternative
approximation scheme is given in \cite{BSchill}.


\section{Approximation}\label{Est}

In this section we apply an approximation scheme recently devised in \cite{Szt2010}. We have
\begin{eqnarray*}
  \Aapprox{\varepsilon} \varphi (x)  
  &  =  &  \int \left(\varphi(y)-\varphi(x)\right)f_{\varepsilon}(x,y)\,dy + 
           (\supmasyL-b_\varepsilon(x))\int (\varphi(y)-\varphi(x))\delta_x(dy) \\
  &  =  &  \int (\varphi(y)-\varphi(x))\tilde{\nu}_{\varepsilon}(x,dy)  \\
  &  =  & \Gamma_\varepsilon \varphi(x) - \supmasyL\varphi(x),\quad \varphi\in \Bbound,\,x\in\R,      
\end{eqnarray*}
where
$$
  \tilde{\nu}_\varepsilon(x,dy)=f_{\varepsilon}(x,y)\,dy+(\supmasyL-b_\varepsilon(x))\delta_x(dy),
$$
and
$$
  \Gamma_\varepsilon \varphi(x)=\int \varphi(y)\tilde{\nu}_\varepsilon(x,dy),\quad \varphi\in \Bbound,x\in\R.
$$
This yields that
\begin{equation}\label{eq:exp}
  e^{t\Aapprox{\varepsilon}}\varphi(x)=e^{t(\Gamma_\varepsilon-\supmasyL I)}\varphi(x)=
  e^{-t\supmasyL}e^{t\Gamma_\varepsilon}\varphi(x).
\end{equation}
A consequence of (\ref{eq:exp}) is that we may consider the operator $\Gamma_\varepsilon$ 
and its powers instead of $\Aapprox{\varepsilon}$.
The fact that $\Gamma_\varepsilon$ is positive enables for more precise estimates.

For $n\in\N$ we define
\begin{eqnarray*}
  f_{n+1,\varepsilon}(x,y)
  &  =   &  \int f_{n,\varepsilon}(x,z) f_{\varepsilon}(z,y)\,dz  \\
  &      &  +\, \left(\supmasyL-b_\varepsilon(y)\right)f_{n,\varepsilon}(x,y) 
            \,+\, \left(\supmasyL-b_\varepsilon(x)\right)^{n}f_{\varepsilon}(x,y),
\end{eqnarray*}
where we let $f_{1,\varepsilon}=f_\varepsilon$.
By induction and Fubini--Tonelli theorem we get
\begin{equation}\label{eq:intf_n}
    \int f_{n,\varepsilon}(x,y)\,dy = \supmasyL^n-\left(\supmasyL-b_\varepsilon(x)\right)^n\,,\quad x\in\R,\,n\in\N.
\end{equation}
\noindent
Also, it was proved in \cite[Lemma 3]{Szt2010} that for all $\varepsilon>0$, $x\in\R$, and $n\in\N$
\begin{equation}\label{eq:GammaDensity}
  \Gamma_\varepsilon^n \varphi(x) = \int \varphi(z) f_{n,\varepsilon}(x,z)\,dz 
  + \left(\supmasyL-b_\varepsilon(x)\right)^{n}\varphi(x),
\end{equation}  
whenever $\varphi\in \Bbound$.

The next lemma is crucial for our further investigation. The significance 
of the inequalities below is that before the expressions on the right hand side we obtain precisely the constants equal to one. 

\begin{lemma}\label{lm:EintL}
We have the following.
\begin{itemize}
\item[(1)] If (A.1), (A.2) and (A.4) hold then there is a constant $c_7=c_7(\phi,\alpha,d)$ and the number $\kappa \in (0,1)$ such that 
$$
\int_{B(y,\kappa|y-x|)} |z-x|^{-\alpha-d} f_{\varepsilon} (y,z) dz \leq \left(b_{\varepsilon}(y) + c_7\right)|y-x|^{-\alpha-d},
$$
for every $\varepsilon \in (0,1)$ and for every $x,y \in \R$.
\item[(2)] If (A.1) and (A.2) hold then there is a constant $c_8=c_8(\phi,\alpha,d)$ such that
$$
\int_{B(y,1)} \frac{\phi(|z-x|)}{|z-x|^{\alpha+d}} f_{\varepsilon} (y,z) dz \leq \left(b_{\varepsilon}(y) + c_8\right)\frac{\phi(|y-x|)}{|y-x|^{\alpha+d}},
$$
for every $\varepsilon \in (0,1)$ and for every $|x-y|>2$.
\end{itemize}
\end{lemma}

\dowod 
First we prove the statement (1). We have 
\begin{align*}
\int_{B(y,\kappa|y-x|)} & |z-x|^{-\alpha-d}f_{\varepsilon}(y,z)\,dz \\ & = \int_{B(y,\kappa|y-x|)} \left[|z-x|^{-\alpha-d}-|y-x|^{-\alpha-d}\right] f_{\varepsilon} (y,z) dz \\
& + |y-x|^{-\alpha-d}  \int_{B(y,\kappa|y-x|)} f_{\varepsilon} (y,z)dz.
\end{align*}
We only need to estimate the first integral on the right hand side of the above equality. Denote $\theta(z):= |z-x|^{-\alpha-d}$, $|z-x|>0$.
\begin{equation}\label{eq:1derphi}
  \partial_j \theta(z)=(\alpha+d)|z-x|^{-\alpha-d-2}(x_j-z_j),
\end{equation}
and
\begin{displaymath}
  \partial_{j,k} \theta(z) = (\alpha+d)|z-x|^{-\alpha-d-2}\left[(\alpha+d+2)\frac{(x_j-z_j)(x_k-z_k)}{|x-z|^{2}}-\delta_{jk}\right].
\end{displaymath}
This yields
\begin{equation}\label{eq:2derphi}
  \sup_{\substack{ z\in B(y,\kappa|y-x|),\\ j,k\in\{1,\dots,d\} }}  |\partial_{j,k} \theta(z)| \leq (\alpha+d)(\alpha+d+3)(1-\kappa)^{-\alpha-d-2}|y-x|^{-\alpha-d-2}, 
\end{equation}
for every $\kappa \in (0,1)$.
Using the Taylor expansion for $\theta$, \eqref{eq:1derphi} and \eqref{eq:2derphi}, (A.1), (A.2) and (A.4), we get
\begin{eqnarray*}
  &      & \left| \int_{B(y,\kappa|y-x|)}\left[ |z-x|^{-\alpha-d} - |y-x|^{-\alpha-d}\right]f_{\varepsilon}(y,z) \,dz\right| \\
  &   =  & \left| \int_{B(0,\kappa|y-x|)} 
                  \left(\theta(y+h)-\theta(y)\right) f_{\varepsilon}(y,y+h)\,dh \right| \\
  & \leq & \left| \int_{B(0,\kappa|y-x|)} 
                  \left(\theta(y+h)-\theta(y)-\scalp{\nabla\theta (y)}{h}\right)
                   f_{\varepsilon}(y,y+h)\,dh \right| \\
  &   +  & \left| \int_{B(0,\kappa|y-x|)} \scalp{\nabla\theta (y)}{h}\,
                    \frac{f_{\varepsilon}(y,y+h)-f_{\varepsilon}(y,y-h)}{2}\,dh \right| \\
  & \leq & C |y-x|^{-\alpha-d}|y-x|^{-\alpha}{\kappa}^{1-\alpha}(\kappa(1-\kappa)^{-\alpha-d-2}+1) \\ 
  & \leq & |y-x|^{-\alpha-d}\left(\indyk{s<\varepsilon_0}(|y-x|) \int_{|z-y|>\kappa|y-x|} f(y,z)\,dz + c_7 \indyk{s \geq \varepsilon_0}(|y-x|)\right),
\end{eqnarray*}
for sufficiently small $\kappa \in (0,1)$. This ends the proof of (1).

We now show the statement (2). Let $|x-y|>2$. Similarly as before we have 
\begin{align*}
\int_{B(y,1)} \frac{\phi(|z-x|)}{|z-x|^{\alpha+d}}f_{\varepsilon}(y,z)\,dz & = \int_{B(y,1)} \left[\frac{\phi(|z-x|)}{|z-x|^{\alpha+d}}-\frac{\phi(|y-x|)}{|y-x|^{\alpha+d}}\right] f_{\varepsilon} (y,z) dz \\
& + \frac{\phi(|y-x|)}{|y-x|^{\alpha+d}}  \int_{B(y,1)} f_{\varepsilon} (y,z)dz.
\end{align*}
Observe that it is enough to estimate the first integral on the right hand side of the above-displayed equality. Denote $\eta(z):= \phi(|z-x|)|z-x|^{-\alpha-d}$. 
Clearly, by (A.1) (a)-(b), we have 
\begin{equation}\label{eq:3derphi}
  \max \left(\sup_{\substack{ z\in B(y,1),\\ j\in\{1,\dots,d\} }}  |\partial_{j} \eta(z)|, \sup_{\substack{ z\in B(y,1),\\ j,k\in\{1,\dots,d\} }}  |\partial_{j,k} \eta(z)| \right)\leq C \eta(y).
\end{equation}
Using the Taylor expansion for $\eta$, (\ref{eq:3derphi}), (A.1) and (A.2), we obtain
\begin{eqnarray*}
  &      & \left| \int_{B(y,1)}\left[\frac{\phi(|z-x|)}{|z-x|^{\alpha+d}}-\frac{\phi(|y-x|)}{|y-x|^{\alpha+d}}\right] f_{\varepsilon} (y,z) \,dz\right| \\
  &   =  & \left| \int_{B(0,1)} 
                  \left(\eta(y+h)-\eta(y)\right) f_{\varepsilon}(y,y+h)\,dh \right| \\
  & \leq & \left| \int_{B(0,1)} 
                  \left(\eta(y+h)-\eta(y)-\scalp{\nabla\eta (y)}{h}\right)
                    f_{\varepsilon}(y,y+h)\,dh \right| \\
  &   +  & \left| \int_{B(0,1)} \scalp{\nabla\eta (y)}{h}\,
                    \frac{f_{\varepsilon}(y,y+h)-f_{\varepsilon}(y,y-h)}{2}\,dh \right| \\
  & \leq & c_8 \eta(y).
\end{eqnarray*}
which ends the proof.
\qed

We now obtain estimates of $f_{n,\varepsilon}(x,y)$. Our argument in the proof of the following lemma shows significance of assumptions on the dominating function $\phi$.  

\begin{lemma}\label{lm:Estimate1}
  If (A.1) -- (A.4) hold then:
\begin{itemize}
\item[(1)] there exists a constant $c_9=c_9(\phi,\alpha,d)$ such that
$$
    f_{n,\varepsilon}(x,y) \leq c_9 n\left(\supmasyL + c_7 \right)^{n-1}|y-x|^{-\alpha-d},
$$
for every $x,y\in\R$, $\varepsilon \in (0,1)$, $n \in \N$,

\item[(2)] there exist the constants $c_{10}=c_{10}(\phi,\alpha,d)$ and $c_{11}=c_{11}(\phi,\alpha,d)$ such that
$$
    f_{n,\varepsilon}(x,y) \leq c_{10} n\left(\supmasyL + c_{11}\right)^{n-1}\frac{\phi(|y-x|)}{|y-x|^{\alpha+d}},
$$
for every $x,y\in\R$, $\varepsilon \in (0,1)$, $n \in \N$.
\end{itemize}
\end{lemma} 

\dowod We use induction. Clearly, for $n=1$ both inequalities hold with constants $c_9 = M \,$, $c_{10} = M$ (and an arbitrary positive $c_{11}$), respectively. Consider first the inequality in (1). We will prove that it holds with constant $c_9 = M \kappa^{-\alpha-d} \,$, where $\kappa \in (0,1)$ is the number from previous lemma. 

Let 
$$
\int f_{n,\varepsilon}(x,z) f_{\varepsilon}(z,y) dz = \int_{B(y,\kappa|y-x|)^c} + \int_{B(y,\kappa|y-x|)} = I + II.
$$
By (A.1) (a) and \eqref{eq:intf_n}, we have 
\begin{align*}
I \leq \kappa^{-\alpha-d} M \,&  |y-x|^{-\alpha-d}  \int f_{n,\varepsilon}(x,z)dz \\ & = \kappa^{-\alpha-d}M \, |y-x|^{-\alpha-d} \left[\bar{b}^n_{\varepsilon}-(\bar{b}_{\varepsilon}-b_{\varepsilon}(x))^n\right].
\end{align*}
By symmetry of $f$ (see (A.3)), induction and Lemma \ref{lm:EintL} (1), we also have
\begin{align*}
II \leq c_9 n (\supmasyL + c_7)^{n-1} & \int_{B(y,\kappa|y-x|)} |x-z|^{-\alpha-d}f_{\varepsilon}(y,z) dz \\ & \leq  c_9 n (\supmasyL + c_7)^{n-1} (b_{\varepsilon}(y)+c_7)|x-y|^{-\alpha-d}.
\end{align*}
We get
\begin{eqnarray*}
  f_{n+1,\varepsilon}(x,y)
  &  =   & I+II + \left(\supmasyL-b_\varepsilon(y)\right)f_{n,\varepsilon}(x,y) +
          \left(\supmasyL-b_\varepsilon(x)\right)^{n}f_{\varepsilon}(x,y) \\
  & \leq &  M \kappa^{-\alpha-d}\,  
           \left[\supmasyL^n-\left(\supmasyL-b_\varepsilon(x)\right)^n\right]|y-x|^{-\alpha-d}\\
  &      & + c_9 n (\supmasyL + c_7)^{n-1} (b_{\varepsilon}(y)+c_7)|x-y|^{-\alpha-d}   \\
  &      & + \left(\supmasyL-b_\varepsilon(y)\right)c_9 n (\supmasyL + c_7)^{n-1}
           |x-y|^{-\alpha-d}   \\
  &      & + \left(\supmasyL-b_\varepsilon(x)\right)^{n}M|x-y|^{-\alpha-d} \\
  &   \leq  & c_9 (n+1) (\supmasyL + c_7)^{n}|x-y|^{-\alpha-d},
\end{eqnarray*}
which ends the proof of part (1).

We now complete the proof of the inequality in (2). We will prove that it holds with constants $c_{10} = c_1 \max(c_9, 2^{\alpha+d}M)$ and $c_{11} = \max(c_7, c_8 + Mc_3)$. When $|x-y| \leq 2$, then it directly follows from the part (1) and (A.1)(a). Assume now that $|x-y| > 2$.  
We have
$$
\int f_{n,\varepsilon}(x,z) f_{\varepsilon}(z,y) dz = \int_{B(x,1)} + \int_{B(x,1)^c} = I + II.
$$
By (A.1) (a) and \eqref{eq:intf_n}, we get 
\begin{align*}
I \leq 2^{\alpha+d} M c_1 \frac{\phi(|x-y|)}{|y-x|^{\alpha+d}} & \int f_{n,\varepsilon}(x,z)dz \\ & = 2^{\alpha+d}M c_1 \frac{\phi(|x-y|)}{|y-x|^{\alpha+d}} \left[\bar{b}^n_{\varepsilon}-(\bar{b}_{\varepsilon}-b_{\varepsilon}(x))^n\right].
\end{align*}
By symmetry of $f$ (see (A.3)), induction, Lemma \ref{lm:EintL} (2) and (A.1) (c), we also have
\begin{align*}
II \leq c_{10} n (\supmasyL + c_{11})^{n-1} & \int_{B(x,1)^c} \frac{\phi(|x-z|)}{|x-z|^{\alpha+d}}f_{\varepsilon}(y,z) dz \\ & \leq  c_{10} n (\supmasyL + c_{11})^{n-1} (b_{\varepsilon}(y)+c_8+Mc_3)\frac{\phi(|x-y|)}{|x-y|^{\alpha+d}}.
\end{align*}
We get
\begin{eqnarray*}
  f_{n+1,\varepsilon}(x,y)
  &  =   & I+II + \left(\supmasyL-b_\varepsilon(y)\right)f_{n,\varepsilon}(x,y) +
          \left(\supmasyL-b_\varepsilon(x)\right)^{n}f_{\varepsilon}(x,y) \\
  & \leq &  2^{\alpha+d}M c_1  
           \left[\supmasyL^n-\left(\supmasyL-b_\varepsilon(x)\right)^n\right] \frac{\phi(|y-x|)}{|y-x|^{\alpha+d}}\\
  &      & + c_{10} n (\supmasyL + c_{11})^{n-1} (b_{\varepsilon}(y)+c_8+Mc_3)\frac{\phi(|x-y|)}{|x-y|^{\alpha+d}}   \\
  &      & + \left(\supmasyL-b_\varepsilon(y)\right)c_{10} n (\supmasyL + c_{11})^{n-1}
           \frac{\phi(|x-y|)}{|x-y|^{\alpha+d}}   \\
  &      & + \left(\supmasyL-b_\varepsilon(x)\right)^{n}M\frac{\phi(|x-y|)}{|x-y|^{\alpha+d}} \\
  &   \leq  & c_{10} (n+1) (\supmasyL + c_{11})^{n}\frac{\phi(|x-y|)}{|x-y|^{\alpha+d}}.
\end{eqnarray*}
\qed

\begin{lemma}\label{lm:Estimate2} Assume (A.1), (A.3) and (A.4). Then
there exists $c_{12}=c_{12}(\phi,\alpha,d)$ such that
\begin{equation}\label{eq:Estimate2}
    f_{n,\varepsilon}(x,y)\leq c_{12} \supmasyL^{d/\alpha}\left(\supmasyL^n-\left(\supmasyL-b_\varepsilon(x)\right)^n\right),
    \quad x,y\in\R,\,\varepsilon\in(0, \varepsilon_0 ),\,n\in\N.
  \end{equation}  
\end{lemma}

\dowod For $n=1$ by (A.1) and (A.4) we have 
$$
  f_{\varepsilon}(x,y)\leq \frac{M \,}{\varepsilon^{\alpha+d}} \leq 
  M \,\left( \frac{b_\varepsilon(x)}{c_6} \right)^{(\alpha+d)/\alpha}
  \leq
  M \, 
  \left(\frac{b_\varepsilon(x)}{c_6} \right)
  \left(\frac{\supmasyL}{c_6}\right)^{d/\alpha}  ,
$$ 
and so (\ref{eq:Estimate2}) holds with $c_{12}=M \,{c_6}^{-d/\alpha-1}$.
Let (\ref{eq:Estimate2}) holds for some $n\in\N$ with $c_{12}=M \,{c_6}^{-d/\alpha-1}$. 
By induction and the symmetry of 
$f_{\varepsilon}$ we get
\begin{eqnarray*}
  f_{n+1,\varepsilon}(x,y)
  & \leq & c_{12} \supmasyL^{d/\alpha}\left(\supmasyL^n-\left(\supmasyL-b_\varepsilon(x)\right)^n\right)\left(\int f_{\varepsilon}(y,z)\,dz + \supmasyL-b_\varepsilon(y)\right)\\
  &      & + \left(\supmasyL-b_\varepsilon(x)\right)^n c_{12} \supmasyL^{d/\alpha}b_\varepsilon(x) \\
  &   =  & c_{12} (\supmasyL)^{d/\alpha}\left(\supmasyL^{n+1}-\left(\supmasyL-b_\varepsilon(x)\right)^{n+1}\right).
\end{eqnarray*}
\qed

In the following lemma we will need some additional notation. For a function $g$ we denote: $b^g_{\varepsilon}(x) := \int_{|y-x|>\varepsilon} g(|y-x|) f_{\varepsilon}(x,y) dy$ and $\bar{b}^g_{\varepsilon} = \sup_{x \in \R} b^g_{\varepsilon}(x)$.
 We note that it follows from (A.1) that 
$$
  \bar{b}_\varepsilon^{\frac{1}{\phi}} \leq c_{13} \varepsilon^{-\alpha}.
$$

\begin{lemma}\label{lm:Estimate3} If (A.1), (A.3) and (A.4) are satisfied then there exist $c_{14}=c_{14}(\phi,\alpha,d)$ and $c_{15}=c_{15}(\phi,\alpha,d)$ such that
\begin{equation}\label{eq:Estimate3}
    f_{n,\varepsilon}(x,y)\leq c_{14} \left(\supmasyL+c_{15}\right)^{n+d/\alpha}
    n^{-d/\alpha},
    \quad x,y\in\R,\,\varepsilon\in(0,\varepsilon_0 \wedge 1),\,n\in\N.
  \end{equation}  
\end{lemma}
\dowod 
We may choose $n_0\in\N$ such that 
\begin{equation}\label{eq:defn_0}
  (1-c_6/c_5)^n (n+1)^{d/\alpha}<\frac{1}{n+1}
\end{equation} 
for every $n\geq n_0$. For $n\leq n_0$ by Lemma \ref{lm:Estimate2} 
we have 
$$
  f_{n,\varepsilon}(x,y)\leq c_{12} \supmasyL^{d/\alpha}\supmasyL^n \leq c_{12} \supmasyL^{n+d/\alpha}
    n^{-d/\alpha} n_0^{d/\alpha},
$$
which yields
the inequality (\ref{eq:Estimate3}) with $c_{14}=c_{12}n_0^{d/\alpha}$ in this case. For $n\geq n_0$ we use induction.
We assume
that (\ref{eq:Estimate3}) holds for some $n\geq n_0$ with
$c_{14}=\max(c_{12}n_0^{d/\alpha},M \, \eta^{-\alpha-d}{c_6}^{-1-d/\alpha})$ and $c_{15}= \bar{b}^{\frac{1}{\phi}-1}_1$, where
$$
p=\frac{d\,2^{\max(d/\alpha,1)-1}}{\alpha},\quad \mbox{and}  \quad
  \eta=\left(\frac{c^2_4/(c_{13}(c_5+c_{15}))}{2+2p}\right)^{\frac{1}{\alpha}}.
$$
 We have
$$
  \int  f_{n,\varepsilon}(x,z)f_{\varepsilon}(z,y)\,dz = \int_{B(y,\eta \varepsilon(n+1)^{1/\alpha})^c } + \int_{B(y,\eta \varepsilon(n+1)^{1/\alpha})}
  = I+II.
$$
By (A.1), (A.4) and (\ref{eq:intf_n}) we get
\begin{eqnarray*}
  I
  &   =  & \int_{B(y,\eta \varepsilon(n+1)^{1/\alpha})^c} f_{n,\varepsilon}(x,z) f_{\varepsilon}(z,y)\,dz \\
  & \leq & M \,\int_{B(y,\eta \varepsilon(n+1)^{1/\alpha})^c} f_{n,\varepsilon}(x,z)|y-z|^{-\alpha-d}\,dz \\
  & \leq & M \,\eta^{-\alpha-d}\varepsilon^{-\alpha-d} (n+1)^{-1-d/\alpha}\int f_{n,\varepsilon}(x,z)\,dz \\
  & \leq & M \,\eta^{-\alpha-d}{c_6}^{-1-d/\alpha}\supmasyL^{1+d/\alpha}(n+1)^{-1-d/\alpha}
           \left[\supmasyL^n-\left(\supmasyL-b_\varepsilon(x)\right)^n\right].
\end{eqnarray*}
By induction, the symmetry of $f_{\varepsilon}$ and (A.4) we obtain
\begin{eqnarray*}
  II
  &   =  & \int_{B(y,\eta \varepsilon(n+1)^{1/\alpha})} f_{n,\varepsilon}(x,z)f_{\varepsilon}(z,y)\,dz \\
  & \leq & c_{14} \left(\supmasyL+c_{15}\right)^{n+d/\alpha}n^{-d/\alpha} 
           \int_{B(y,\eta \varepsilon(n+1)^{1/\alpha})} \frac{f_{\varepsilon}(y,z)}{\phi(|z-y|)} \,dz \\
  &   =  & c_{14} \left(\supmasyL+c_{15}\right)^{n+d/\alpha}n^{-d/\alpha} 
          \left(b^{\frac{1}{\phi}}_\varepsilon(y)-b^{\frac{1}{\phi}}_{\eta \varepsilon (n+1)^{1/\alpha}}(y)\right)\\
  & \leq & c_{14} \left(\supmasyL+c_{15}\right)^{n+d/\alpha}n^{-d/\alpha} 
          b^{\frac{1}{\phi}}_\varepsilon(y)\left(1-\frac{c_4\eta^{-\alpha}}{c_{13}(n+1)}\right).
\end{eqnarray*}
By (\ref{eq:defn_0}) we also have 
\begin{equation}\label{eq:hilfe}
  \left(1-\frac{b_\varepsilon(x)}{\supmasyL}\right)^n (n+1)^{d/\alpha} \leq (1-c_6/c_5)^n (n+1)^{d/\alpha} \leq
  \frac{1}{n+1}.
\end{equation}
Using the fact that $\phi(a)=1$ for $a\in [0,1]$ and 
$b^{\frac{1}{\phi}}_\varepsilon(y)-b_\varepsilon(y) = b^{\frac{1}{\phi}-1}_1(y)\leq c_{15}$ we get
\begin{eqnarray*}
  f_{n+1,\varepsilon}(x,y)
  &   =  & I+II + \left(\supmasyL-b_\varepsilon(y)\right)f_{n,\varepsilon}(x,y) +
          \left(\supmasyL-b_\varepsilon(x)\right)^{n}f_{\varepsilon}(x,y) \\
  & \leq & c_{14} \supmasyL^{1+d/\alpha}(n+1)^{-1-d/\alpha}\left[\supmasyL^n-\left(\supmasyL-b_\varepsilon(x)\right)^n\right]\\
  &      & + c_{14} \left(\supmasyL+c_{15}\right)^{n+d/\alpha}n^{-d/\alpha} b^{\frac{1}{\phi}}_\varepsilon(y) 
           \left(1-\frac{c_4\eta^{-\alpha}}{c_{13}(n+1)}\right)\\
  &      & + c_{14} \left(\supmasyL+c_{15}\right)^{n+d/\alpha}n^{-d/\alpha} \left(\supmasyL-b_\varepsilon(y)\right)\\
  &      & + c_{14} \supmasyL^{1+d/\alpha} \left(\supmasyL-b_\varepsilon(x)\right)^{n} \\
  & \leq & c_{14} \left(\supmasyL+c_{15}\right)^{n+1+d/\alpha}(n+1)^{-d/\alpha}
           \left[\frac{1}{n+1}\left(1-\left(1-\frac{b_{\varepsilon}(x)}{\supmasyL}\right)^n\right)\right. \\
  &      & \left. - \frac{b^{\frac{1}{\phi}}_\varepsilon(y)}{\supmasyL+c_{15}}
           \left(1+\frac{1}{n}\right)^{d/\alpha}\frac{c_4\eta^{-\alpha}}{c_{13}(n+1)}
           +\left(1+\frac{1}{n}\right)^{d/\alpha}\right. \\
  &      & \left.+\left(1-\frac{b_\varepsilon(x)}{\supmasyL}\right)^n\left(n+1\right)^{d/\alpha}\right].
\end{eqnarray*}  
By (A.1), (A.4), (\ref{eq:hilfe}) and the following inequality
$$
\frac{b^{\frac{1}{\phi}}_\varepsilon(y)}{\supmasyL+c_{15}} 
\geq \frac{c_4 \varepsilon^{-\alpha}}{c_5 \varepsilon^{-\alpha}+ c_{15}} \geq \frac{c_4 }{c_5+ c_{15}},
$$
the last expression is bounded above by 
\begin{eqnarray*}
  &      & c_{14} \left(\supmasyL+c_{15}\right)^{n+1+d/\alpha}(n+1)^{-d/\alpha} \\
  &      & \times \left[\frac{2}{n+1}+\left(1+\frac{1}{n}\right)^{d/\alpha}
          \left(1-\frac{\eta^{-\alpha}c_4^2/(c_{13}(c_5+c_{15}))}{n+1}\right)\right] \\
  & \leq & c_{14} \left(\supmasyL+c_{15}\right)^{n+1+d/\alpha}(n+1)^{-d/\alpha} \\
  &      & \times \left[\frac{2}{n+1} +
           \left(1+\frac{p}{n}\right)
          \left(1-\frac{\eta^{-\alpha}c_4^2/(c_{13}(c_5+c_{15}))}{n+1}\right)\right] \\
  & \leq &  c_{14} \left(\supmasyL+c_{15}\right)^{n+1+d/\alpha}(n+1)^{-d/\alpha} \\
  &      & \times \left[1 - \frac{1}{n+1}\left( \eta^{-\alpha}c_4^2/(c_{13}(c_5+c_{15}))-2 - 2p\right)  \right],
\end{eqnarray*}
which gives 
$$
 f_{n+1,\varepsilon}(x,y)  \leq  c_{14} \left(\supmasyL+c_{15}\right)^{n+1+d/\alpha}(n+1)^{-d/\alpha}.
$$

\qed

Using the above lemmas we may estimate $\Gamma_\varepsilon^n$ and in consequence also the exponent 
operator $e^{t\Aapprox{\varepsilon}}=e^{-t\supmasyL}e^{t\Gamma_\varepsilon}$.

\begin{lemma}\label{lm:SE1} Assume (A.1) -- (A.4). Then for all $x\in\R$ and all nonnegative 
$\varphi\in\Bbound$ such that $x\notin \supp(\varphi)$ we have
  $$
 e^{t\Aapprox{\varepsilon}}\varphi(x) \leq  c_{10} t\exp(c_{11}t) \int \varphi(y)\frac{\phi(|y-x|)}{|y-x|^{\alpha+d}}\,dy,\quad \varepsilon \in (0,1).
$$
\end{lemma}
\dowod
By \eqref{eq:GammaDensity} and Lemma \ref{lm:Estimate1} for every $\varphi$ such that $x\not\in\supp(\varphi)$ we get
$$
 \Gamma_\varepsilon^n\varphi(x) \leq \int \varphi(y) c_{10} n \left(\supmasyL+c_{11}\right)^{n-1}
 \frac{\phi(|y-x|)}{|y-x|^{\alpha+d}}\,dy,
$$
and
\begin{eqnarray*}
 e^{t\Aapprox{\varepsilon}}\varphi(x) 
 & \leq & c_{10}e^{-t\supmasyL} \sum_{n=1}^\infty 
          \frac{t^n n \left(\supmasyL+c_{11}\right)^{n-1}}{n!}  \int \varphi(y)\frac{\phi(|y-x|)}{|y-x|^{\alpha+d}}\,dy \\
 &   =  & c_{10}e^{-t\supmasyL} t \sum_{n=0}^\infty 
          \frac{t^n\left(\supmasyL+c_{11}\right)^{n}}{n!}  \int \varphi(y)\frac{\phi(|y-x|)}{|y-x|^{\alpha+d}}\,dy \\
 &   =  & c_{10} t\exp\left(c_{11}t\right) \int \varphi(y)\frac{\phi(|y-x|)}{|y-x|^{\alpha+d}}\,dy.
\end{eqnarray*}
\qed

\begin{lemma}\label{lm:SE2} Assume (A.1), (A.3) and (A.4). Then for every nonnegative $\varphi\in \Bbound\cap L_1(\R)$ we have
  $$
 e^{t\Aapprox{\varepsilon}}\varphi(x) \leq c_{14}\exp(c_{15}t) t^{-d/\alpha} \int \varphi(y)\,dy + e^{-tb_\varepsilon(x)}\varphi(x),
  $$
for $x\in\R,\,\varepsilon \in (0,\varepsilon_0\wedge 1),t>0.$
\end{lemma}
\dowod We directly deduce from Lemma \ref{lm:Estimate3} that for every $\varphi\in \Bbound\cap L_1(\R)$ 
$$
 \Gamma_\varepsilon^n\varphi(x) \leq  c_{14} (\supmasyL+c_{15})^{n+d/\alpha}
    n^{-d/\alpha}  \int \varphi(y)\,dy + \left(\supmasyL-b_\varepsilon(x)\right)^{n}\varphi(x),
$$
and, consequently, by \cite[Lemma 9]{Szt2010}, we obtain
\begin{eqnarray*}
 e^{t\Aapprox{\varepsilon}}\varphi(x) 
 & \leq & e^{-t\supmasyL} \left[c_{14} \int \varphi(y)\,dy \sum_{n=1}^\infty 
        \frac{ t^n (\supmasyL+c_{15})^{n+d/\alpha}}{n!n^{d/\alpha}} + 
        e^{t\left(\supmasyL-b_\varepsilon(x)\right)}\varphi(x)\right]         \\
 & \leq & c_{14} \exp(c_{15}t) t^{-d/\alpha} \int \varphi(y)\,dy +  e^{-tb_\varepsilon(x)}\varphi(x). 
 \end{eqnarray*}
 \qed

\dowodof \emph{ of Theorem \ref{th:MainT}}. Let $t>0$, $\varphi\in\Bbound$, and $x\in\R$. Denote: $D = \left\{y \in \R: \phi(|y-x|)|y-x|^{-\alpha-d} < t^{-1-d/\alpha} \right\}$. 
Using Lemma \ref{lm:SE1} for $\indyk{D}\varphi$ and Lemma \ref{lm:SE2} for $\indyk{D^c}\varphi$ we obtain
\begin{eqnarray*}
  e^{t\Aapprox{\varepsilon}}\varphi(x)
  &  =   & e^{t\Aapprox{\varepsilon}}[\indyk{D}\varphi](x) + e^{t\Aapprox{\varepsilon}}[\indyk{D^c}\varphi](x) \\
  & \leq & C_1 e^{C_2t}\left[\int_{D} \varphi(y)\frac{t \phi(|y-x|)}{|y-x|^{\alpha+d}}\, dy +
            \int_{D^c} \varphi(y) t^{-d/\alpha}\,dy \right] \\
  &      & + \,e^{-tb_\varepsilon(x)}\varphi(x) \\
  & \leq & C_1 e^{C_2t} \int \varphi(y) \min\left(t^{-d/\alpha},\frac{t\phi(|y-x|)}{|y-x|^{\alpha+d}}\right)\, dy + e^{-tb_\varepsilon(x)}\varphi(x)
\end{eqnarray*}
\qed

\dowodof \emph{ of Theorem \ref{th:main}}. By Lemma 12 in \cite{Szt2010} we have
$$ 
 \lim_{\varepsilon\to 0}\|\gener\varphi - \Aapprox{\varepsilon}\varphi\|_\infty=0
$$
for every $\varphi\in C_\infty^2(\R)$. 
A closure of $\gener$ is a generator of a semigroup and from the Hille-Yosida theorem it follows that the range of $\lambda-\gener$ is dense in $\Cvanish$ 
and therefore by Theorem 5.2 in \cite{Trotter58} (see also \cite{Hasegawa}) we get
$$
  \lim_{\varepsilon\downarrow 0} \|e^{t\Aapprox{\varepsilon}}\varphi-P_t\varphi\|_\infty =0,
$$
for every $\varphi\in \Cvanish$. By Theorem \ref{th:MainT} this yields
$$
  P_t\varphi(x) \leq C_1 e^{C_2t}  \int 
  \varphi(z) \min\left(t^{-d/\alpha},\frac{t  \phi(|z-x|)}{|z-x|^{\alpha+d}}\right) dz,
$$
for every nonnegative $\varphi\in \Cvanish$.
\qed

\section{Discussion of examples}\label{sec:Examples}

We now prove the condition (A.1) (c) for functions $\phi$ of the form \eqref{eq:specphiform} for restricted set of parameters $\beta$ and $\gamma$. First we recall some well known geometric fact, see e.g. \cite[Lemma 5.3]{KulSiu}.  

\begin{lemma}\label{lm:volumeest} The volume of intersection of two balls $B(x,p+k)$ and $B(y,n-p)$ such that $|y-x|=n \in \N$, $1 \leq p \leq n-1$, $0<k \leq n-p$, is less than $c k^{\frac{d+1}{2}}\left(\min\left\{p+k, n-p\right\}\right)^{\frac{d-1}{2}}$. 
\end{lemma} 

\begin{proposition}
Let the function $\phi$ be of the form \eqref{eq:specphiform}. Then the assumption (A.1) (c) is satisfied if $\beta \in (0,1]$ and $\gamma < d/2+\alpha-1/2$.
\end{proposition}

\dowod
Let $\beta \in (0,1]$ and $\gamma < d/2+\alpha-1/2$. First note that there is an absolute constant $C=C(m,\beta,\gamma)$ such that $\phi(s)s^{-d-\alpha} \leq C \phi(u)u^{-d-\alpha}$ for $|s-u|\leq 1$ whenever $s, u \geq 1$. By this fact, with no loss of generality we may and do consider only the case when $|x-y|=n$ for some even natural number $n \geq 4$. Let
\begin{align*}
\int_{B(x,1)^c \cap B(y,1)^c} & \frac{\phi(|y-z|)}{|y-z|^{\alpha+d}} \frac{\phi(|z-x|)}{|z-x|^{\alpha+d}}\,dz \\ & \leq 2\int_{B(x,1)^c \cap B(y,n-1)^c} + \int_{(B(x,1)^c \cap B(x,n-1)) \cup (B(y,1)^c \cap B(y,n-1))} \\ & = 2I + II.
\end{align*}
We have
$$
I \leq \frac{\phi(n-1)}{(n-1)^{\alpha+d}} \int_{B(0,1)^c} \frac{\phi(|z|)}{|z|^{\alpha+d}}\,dz \leq C \frac{\phi(|y-x|)}{|y-x|^{\alpha+d}}
$$
with some constant $C=C(m,\beta,\gamma,\alpha,d)$. 

To estimate the term II we will need the additional notation. For $1\leq p<n/2$ and $0 \leq k < n-p$ we denote: 
\begin{itemize}
\item $D_p := \left\{z \in \R: n-p-1 \leq |z-y| < n-p, |x-z| < |y-z|\right\}$,
\item $D_{p,k}=D_p \cap \left\{z \in \R: p+k \leq |z-x|< p+k+1\right\}$,
\item $n_p := \max\left\{k \in \N: D_{p,k} \neq \emptyset \right\}$.
\end{itemize}
Clearly, $D_p \subset \bigcup_{k=0}^{n_p} D_{p,k}$ and $D_{p,k} \subset B(x,p+k+1) \cap B(y,n-p)$. We have
\begin{align*}
II & \leq 2^{\alpha+d-\gamma+1} |y-x|^{-\alpha-d+\gamma} \int_{1 \leq |y-z| < n-1, |x-z| < |y-z|} \frac{e^{-m |x-z|^{\beta}} e^{-m |y-z|^{\beta}}}{|x-z|^{\alpha+d-\gamma}} dz \\
& = 2^{\alpha+d-\gamma+1} |y-x|^{-\alpha-d+\gamma} \sum_{p=1}^{n/2-1} \int_{D_p} \frac{e^{-m |x-z|^{\beta}} e^{-m |y-z|^{\beta}}}{|x-z|^{\alpha+d-\gamma}} dz \\
& \leq 2^{\alpha+d-\gamma+1} |y-x|^{-\alpha-d+\gamma} \sum_{p=1}^{n/2-1} \sum_{k=0}^{n_p} \int_{D_{p,k}} \frac{e^{-m |x-z|^{\beta}} e^{-m |y-z|^{\beta}}}{|x-z|^{\alpha+d-\gamma}} dz \\
& \leq 2^{\alpha+d-\gamma+1} |y-x|^{-\alpha-d+\gamma}\sum_{p=1}^{n/2-1} \sum_{k=0}^{n_p}  \frac{e^{-m (p+k)^{\beta}} e^{-m (n-p-1)^{\beta}}}{(p+k)^{\alpha+d-\gamma}} \ |D_{p,k}|.
\end{align*}
Notice that $(n-p)^{\beta}-(n-p-1)^{\beta} \leq \beta$ when $\beta \in (0,1]$. Furthermore, since $p+k \leq p+n_p < n - p$, we also have $k^{\beta} + n^{\beta} \leq (p+k)^{\beta} + (n-p)^{\beta}$. These inequalities and Lemma \ref{lm:volumeest} thus yield
\begin{align*}
II & \leq C \frac{e^{-m n^{\beta}}}{|y-x|^{\alpha+d-\gamma}}\sum_{p=1}^{n/2-1} \sum_{k=0}^{n_p}  e^{-m k^{\beta}}k^{\frac{d+1}{2}} (p+k)^{-\alpha-d+\gamma} \  (p+k)^{\frac{d-1}{2}} \\
& \leq C \frac{e^{-m |y-x|^{\beta}}}{|y-x|^{\alpha+d-\gamma}}\sum_{p=1}^{\infty} p^{-\frac{d+1}{2}-\alpha+\gamma} \  \sum_{k=0}^{\infty}  e^{-m k^{\beta}} k^{\frac{d+1}{2}},
\end{align*}
for some $C=C(m,\beta,\gamma,\alpha,d)$. We conclude by observing that for $\beta > 0$ and $\gamma < d/2+\alpha-1/2$ the last two sums are bounded by constant.
\qed

\begin{remark}
\label{rm:rm} 
\begin{itemize}
\item[(1)] When $\beta>1$, then the condition (c) in assumption (A.1) fails. This can be shown by estimating from below the integral 
$$
\int_{B((x+y)/2,1)} \frac{\phi(|y-z|)}{|y-z|^{\alpha+d}} \frac{\phi(|z-x|)}{|z-x|^{\alpha+d}}\,dz
$$
for $|y-x|$ big enough.
\item[(2)] Also, if $\beta = 1$ and $\gamma = d/2+\alpha-1/2$, then at least for $d=1$ the condition (c) in assumption (A.1) does not hold. In this case we have
$$
\int_{1}^{x-1} e^{-(x-z)} (x-z)^{-1} e^{-z} z^{-1} \,dz = 2 \log(x-1) e^{-x} x^{-1}, \quad  x > 2. 
$$
\end{itemize}
\end{remark}

\end{document}